\documentclass[]{amsart}
\usepackage{amssymb,amscd%,showkeys%
}
\newcommand\commentout[1]{}
\parskip \medskipamount
\newcommand\aDPr{\fraka^D_{\Prin}} % sheaf of principal adeles on S

\newcommand\ad{\operatorname{ad}}

\newcommand\bbeta{{\boldsymbol \beta}}
\newcommand\bcdot{{\raise0.2ex \hbox{\bf.}}}
\newcommand\bdot{{\raise0.4ex \hbox{\bf.}}}

\newcommand\bgamma{{\boldsymbol \gamma}}

\newcommand\bos{{{\mathrm b}{\mathrm o}{\mathrm s}}}
\newcommand\Bos{\operatorname{Bos}}

\newcommand\bsigma{{\boldsymbol \sigma}}
\newcommand\bzeta{{\boldsymbol \zeta}}

\newcommand\Complex{{\mathbb C}}

\newcommand\der{\partial}

\newcommand\fraka{{\mathfrak a}}
\newcommand\free{{{\mathrm{free}}}}

\newcommand\g{{\mathfrak g}}

\newcommand\germanS{{\mathfrak S}}
\newcommand\gh{{{\mathrm{gh}}}}
\newcommand\Gh{\operatorname{Gh}\nolimits}

\newcommand\ghFock{{\mathcal F}^{\gh}}
\newcommand\ghvac{|0\rangle^{\gh}}

\newcommand\gDPr{\g^D_{\Prin}} % sheaf of principal adeles on S

\newcommand\h{{\mathfrak h}}

\newcommand\HOM{\mathop{{\mathcal H}{\mathit o}{\mathit m}}\nolimits}
\newcommand\Hom{\operatorname{Hom}\nolimits}
\newcommand\hor{{\mathrm{hor}}}
\newcommand\hvee{h^\vee}

\newcommand\Ind{\operatorname{Ind}\nolimits}
\newcommand\Integer{{\mathbb Z}}

\newcommand\isoot{\overset\sim\leftarrow}
\newcommand\khat{{\hat k}} % center of aff. Lie alg.

\newcommand\M{{\mathcal M}}
\ifx\mathit\undefined
\newcommand\mathit{} % for old amslatex
\fi
\newcommand\N{{\mathcal N}}
\newcommand\n{{\mathfrak n}}

\newcommand\np{{:}} % normal product
\newcommand\NP{{\genfrac{}{}{0pt}{1}{\circ}{\circ}}} % Normal Product
\renewcommand\O{{\mathcal O}}

\ifx\pd\undefined
\newcommand\pd[2]{\frac{\partial #1}{\partial #2}}
\else
\renewcommand\pd[2]{\frac{\partial #1}{\partial #2}} % for old amslatex
\fi

\newcommand\Prin{{\dot \X}} % principal adele symbol (outside of divisor)

\newcommand\Res{\mathop{\operatorname{Res}}\nolimits}

\newcommand\scr{scr}
\newcommand\sCB{{{\mathcal C}{\mathcal B}}}% Sheaf of Conformal Blocks
\renewcommand\setminus{\smallsetminus}
\newcommand\simeqq{\cong}

\newcommand\tensor{\otimes}

\newcommand\Tr{\operatorname{Tr}\nolimits}
\newcommand\Vac{{{\mathrm V}ac}}

\newcommand\Waki{{\mathrm W}ak}

\newcommand\X{{\mathfrak X}}

\newtheorem{thm}{Theorem}[section]
\newtheorem{lem}[thm]{Lemma}
\newtheorem{prop}[thm]{Proposition}

\theoremstyle{definition}
\newtheorem{defn}[thm]{Definition}

\theoremstyle{remark}
\newtheorem{example}[thm]{Example}

\numberwithin{equation}{section}

\newcommand\thmref[1]{Theorem~\ref{#1}}

\newcommand\secref[1]{Section~\ref{#1}}

\newcommand\lemref[1]{Lemma~\ref{#1}}
\newcommand\defref[1]{Definition~\ref{#1}}

\newcommand\exref[1]{Example~\ref{#1}}

\begin{document}
%%%%%%%%%%%%%%%%%%%%%%%%%%%%%%%%%%%%%%%%%%%%%%%%%%%%%%%%%%%%%%%%%%%%%%%%%%%

\title[WZW Model on elliptic curves]
{Wess-Zumino-Witten model on elliptic curves\\
at the critical level}

%%%%%%%%%%%%%%%%%%%%%%%%%%%%%%%%%%%%%%%%%%%%%%%%%%%%%%%%%%%%%%%%%%%%%%%%%%%

\author{Gen KUROKI} %
\address{Mathematical Institute, Tohoku University, Sendai 980, JAPAN}
\author{Takashi TAKEBE}%
\address{Department of Mathematics, Faculty of Sciences,
Ochanomizu University, Otsuka, Tokyo 112-8610, JAPAN}

\date{Talk given at Mathematical Methods of Regular Dynamics dedicated
to the 150$^{\rm th}$ anniversary of Sophie Kowalevski, University of
Leeds, 12-15 April, 2000}

\ifx\ClassWarning\undefined
\maketitle % for old amslatex
\fi

\begin{abstract}
 We construct a Gaudin type lattice model as the Wess-Zumino-Witten
 model on elliptic curves at the critical level. Bethe eigenvectors are
 obtained by the bosonisation technique.
\end{abstract}

\ifx\ClassWarning\undefined\else
\maketitle
\fi

%%%%%%%%%%%%%%%%%%%%%%%%%%%%%%%%%%%%%%%%%%%%%%%%%%%%%%%%%%%%%%%%%%%%%%%%%%%
% for debug
%%%%%%%%%%%%%%%%%%%%%%%%%%%%%%%%%%%%%%%%%%%%%%%%%%%%%%%%%%%%%%%%%%%%%%%%%%%
%\centerline{\VERSION.}
%\centerline{\bf \fbox{Not for distribution!}}

%%%%%%%%%%%%%%%%%%%%%%%%%%%%%%%%%%%%%%%%%%%%%%%%%%%%%%%%%%%%%%%%%%%%%%%%%%%
% for debug
%%%%%%%%%%%%%%%%%%%%%%%%%%%%%%%%%%%%%%%%%%%%%%%%%%%%%%%%%%%%%%%%%%%%%%%%%%%
\setcounter{section}{-1}
%%%%%%%%%%%%%%%%%%%%%%%%%%%%%%%%%%%%%%%%%%%%%%%%%%%%%%%%%%%%%%%%%%%%%%%%%%%
\section{Introduction}
\label{sec:intro}

The goal of this article is to construct a lattice model which is a
variant of the Gaudin model with the help of the Wess-Zumino-Witten
(WZW) model on elliptic curves at the critical level and to find its
eigenvectors by means of the bosonisation of the WZW model.

As is well known, correlation functions of the WZW model on the Riemann
sphere satisfy the Knizhnik-Zamolodchikov equations when the level $k$
of the model is not critical, i.e., $k\neq -h^\vee$, where $h^\vee$ is
the dual Coxeter number of the simple Lie algebra $\g$ which describes
symmetry of the model. In constrast to this case, when the level is
critical, $k = - h^\vee$, there can be no longer such an equation since
the Sugawara construction of the energy-momentum tensor breaks down.

In \cite{f-f-r:94} Feigin, Frenkel and Reshetikhin found an
interpretation of such case as a lattice model called the (rational)
Gaudin model \cite{gau:73}, \cite{gau:76}, \cite{gau:83}, a
quasi-classical limit of the totally isotropic spin chain model (the XXX
model). (See \cite{skl:87}.) It is also shown in \cite{f-f-r:94} that
the free field realisation of the WZW model provides a new
diagonalisation method which recovers the results of the Bethe Ansatz
method.

In this paper we apply this story to the WZW model on elliptic
curves. The state space of the lattice model thus obtained is a space of
functions over the Cartan subalgebra. The transfer matrix (the
generating function of the Hamiltonians) is a quasi-classical limit of
the IRF type lattice model. Therefore we name the model a ``face type
Gaudin model''. The bosonisation technique is also applied to find a
Bethe Ansatz type eigenvectors.

Felder and Varchenko \cite{fel-var:95}, \cite{fel-var:97} studied this
system and its Bethe Ansatz, which arose from the stationary phase
method of the integral representation of solutions of the KZB equations. 
Their results are used in the analysis of a spin chain model with
elliptic exchanges in \cite{ino:99}. The same kind of system (the
Gaudin-Calogero model) has been also studied by Enriquez, Feigin and
Rubtsov \cite{enr-rub:96}, \cite{e-f-r:98} who started from the
quantisation of the Hitchin system on elliptic curves.

This paper is organised as follows. 
In \secref{sec:result} we state our main
results. \secref{sec:wzw-critical} explains how to derive our transfer
matrix from the WZW model on elliptic curves at the critical level. In
order to construct the Bethe eigenvectors we make use of the Wakimoto
modules with the critical level, which we recall in
\secref{sec:wakimoto-critical}. The last section,
\secref{sec:free-field,ba}, shows that the free field theory gives the
eigenvector of the transfer matrix in the form of the Bethe vector.

Details shall be published in a forthcoming paper. Mathematically
delicate conditions like finiteness of modules are not specified unless
they are essential.

%%%%%%%%%%%%%%%%%%%%%%%%%%%%%%%%%%%%%%%%%%%%%%%%%%%%%%%%%%%%%%%%%%%%%%
\section{Main results}
\label{sec:result}

First we fix notations. Throughout this paper $\tau$ is a complex number
with positive imaginary part. The elliptic curve with modulus $\tau$ is
denoted by $X = X_\tau = \Complex/\Integer + \Integer\tau$.

Let $\g$ be a finite dimensional simple Lie algebra of rank $l$, $\h$
be its Cartan subalgebra and
\begin{equation}
    \g = \h \oplus \bigoplus_{\alpha \in \Delta} \g_\alpha
\label{root-decomp}
\end{equation}
be the root space decomposition, where $\Delta$ is the set of roots. We
use the Cartan-Killing form normalised as follows:
\begin{equation}
    (A \,|\, B) :=
    \frac{1}{2 h^\vee} \Tr_{\g}(\ad A \ad B) \in \Complex
    \text{ for $A,B\in\g$},
\label{def:killing}
\end{equation}
We identify $\h$ and its dual space $\h^\ast$ through this inner
product.  We fix the simple roots $\{\alpha_1,\ldots,\alpha_l\}$,
Chevalley generators $\{H_i, E_i, F_i\}_{i=1,\dots,l}$ and a basis
$e_\alpha$ of $\g_\alpha$, such that $e_{\alpha_i} = E_i$ for
$i=1,\ldots,l$ and
$(e_\alpha|e_{-\alpha'})=\delta_{\alpha,\alpha'}$. The set of positive
and negative roots are denoted by $\Delta_+ = \{\beta_1, \dots, \beta_s
\}$ and $\Delta_-$ respectively. We fix an orthonormal basis
$\{h_r\}_{r=1,\dots,l}$ of $\h$ and the coordinate system of $\h$,
$(\xi_1,\dots,\xi_l) \in \Complex^l$, associated to it.

Let $V_i$ ($i=1,\dots,N$) be finite dimensional irreducible
representations of $\g$ with the highest weight $\lambda_i$ and $V(0)$
be the $0$-weight space of $V = V_1 \tensor \cdots \tensor V_N$. We
denote the action of $\g$ on the $i$-th factor $V_i$ by $\rho_i$ as
usual. The dual (right) action of $\g$ on $\Phi \in V^\ast =
\Hom_\Complex(V,\Complex)$ is denoted as
\begin{equation}
    \rho_i^\ast(A) \Phi(v) := \Phi( \rho_i(A)v),
\label{def:rho-ast}
\end{equation}
and the $0$-weight space of $V^\ast$ by $V^\ast(0)$.

We define a differential operator $\hat\tau(u)$ on $V^\ast(0)$-valued
functions on
\begin{equation}
    S=\{ H \in \h \mid 
         \alpha(H) \not\in \Integer \text{ for all } \alpha \in \Delta\}
\label{def:S}
\end{equation}
with fixed complex parameters $z_i$ ($i=1,\dots,N$) and a spectral
parameter $u$ as follows:
\begin{equation}
\begin{split}
    & \hat\tau(u) :=
    \frac{1}{2} \sum_{r=1}^l \nabla_{\xi_r,u}^2
    + \\ &+
    \frac{1}{2} \sum_{i,j=1}^N \sum_{\alpha\in\Delta}
    w_{\alpha(H)}(z_i-u)\, w_{-\alpha(H)}(z_j-u) 
    \rho_j^\ast(e_{-\alpha}) \rho_i^\ast(e_{\alpha}),
\end{split}
\label{def:tau(u)}
\end{equation}
where
\begin{equation}
    \nabla_{\xi_r,u} := 
    \frac{\der}{\der \xi_r} 
    - \sum_{i=1}^N \zeta_{11}(z_i-u) \rho^\ast_i(h_r).
\label{def:nabla}
\end{equation}
The quasi-periodic functions $w_c(z)$ and $\zeta_{11}(z)$ are defined by
\eqref{def:w,zeta}.

Since this operator can be interpreted as a trace of square of the
dynamical (or modified) classical $r$-matrix which is a classical limit
of the IRF type lattice models (cf.\ \cite{fel:94}, \cite{fel:95},
\cite{fel-sch:99}), we call $\hat\tau(u)$ the {\em transfer matrix of
the face type Gaudin model}.

\begin{thm}
\label{thm:commutativity}
 The operator $\hat\tau(u)$ commutes with itself:
\begin{equation}
    [\hat\tau(u), \hat\tau(u')] = 0.
\label{[tau,tau]=0}
\end{equation}
\end{thm}

This can be checked by direct computation, but we shall show in
\secref{sec:wzw-critical} that it can be proved with the help of the WZW
model at the critical level.

Taking \thmref{thm:commutativity} into account, we can pose a question
of simultaneous diagonalisation of $\hat\tau(u)$. Our main result is the
following Bethe Ansatz solution of this problem. Assume that $V_i$ is
the dual Verma module $M^\ast_{\lambda_i}$ of $\g$ and that the sum of
the weights $\lambda_i$ belongs to the positive root lattice:
\begin{equation}
    \sum_{i=1}^N \lambda_i = \sum_{j=1}^M \alpha_{i(j)},
\label{charge-preserving}
\end{equation}
for a sequence $\{\alpha_{i(j)}\}_{j=1,\dots,M}$ of simple roots. The
symbol $\jmath(v)$ for $v\in M^\ast_\lambda$ is the canonical paring of
$v$ with the highest weight vector of the Verma module $M_\lambda$.

\begin{thm}
\label{thm:BA}
 If there are complex numbers $t_j$ ($j=1,\dots,M$) satisfying the {\em
 Bethe Ansatz equation},
\begin{equation}
    \sum_{i=1}^N (\alpha_{i(j)}\mid \lambda_i) \zeta_{11}(t_j-z_i)
    =
    \sum_{j'\neq j} 
    (\alpha_{i(j)}\mid \alpha_{i(j')}) \zeta_{11}(t_j-t_{j'}),
\label{BA-eq}
\end{equation}
for any $j = 1,\dots, M$, then 
\begin{equation}
    \Psi(H;v) :=
    \sum_{\{I_j\}} \prod_{a=1}^N 
    \langle I_a; v_a; z_a, t_j (j\in I_a) \rangle
\label{BA-vector}
\end{equation}
 is an eigenvector of $\hat\tau(u)$. Here $\{I_j\}$ is a partition of the
 set $\{1,\dots,M\}$ into $N$ sets, $I_1 \sqcup \cdots \sqcup I_N =
 \{1,\dots,M\}$, and the symbol $\langle I; v; t_j (j\in I)\rangle$ is
 defined as follows: 
\begin{multline}
    \langle I; v; z, t_j (j\in I) \rangle
    :=
    \sum_{\sigma \in \germanS} 
    w_{\alpha_{i(\sigma(1))}} (t_{\sigma(1)} - t_{\sigma(2)})\,
    w_{\alpha_{i(\sigma(1))} + \alpha_{i(\sigma(2))}} 
    (t_{\sigma(2)} - t_{\sigma(3)}) \times \cdots
\\
    \cdots \times
    w_{\alpha_{i(\sigma(1))}+\alpha_{i(\sigma(2))}+\cdots
      +\alpha_{i(\sigma(m))}} (t_{\sigma(m)} - z)\,
    \jmath( E_{i(\sigma(m))} \cdots E_{i(\sigma(1))} v)
\end{multline}
 if $I = \{1,\dots, m\}$. (We write $w_\alpha$ instead of $w_{\alpha(H)}$ 
 for short.) The eigenvalue of $\Psi(H;v)$ is
\begin{equation}
    \tau_\Psi(u) :=
    \frac{1}{2} \sum_{r=1}^l
    \bzeta(h_r;z,t;u)^2
    + \frac{\der}{\der u} \bzeta(\rho;z,t;u).
\label{BA-eigenval}
\end{equation}
Here $\bzeta$ is defined by
\begin{equation}
    \bzeta(h;z,t;u) :=
    \sum_{i=1}^N \lambda_i(h) \zeta_{11}(z_i - u) 
    -
    \sum_{j=1}^M \alpha_{i(j)}(h) \zeta_{11}(t_j - u),
\label{def:bzeta}
\end{equation}
and $\rho$ is the half sum of the positive roots of $\g$.
\end{thm}

We shall show how to prove this theorem by bosonisation in
\secref{sec:free-field,ba}.

%%%%%%%%%%%%%%%%%%%%%%%%%%%%%%%%%%%%%%%%%%%%%%%%%%%%%%%%%%%%%%%%%%%%%%
\section{WZW model on elliptic curves at the critical level}
\label{sec:wzw-critical}

The idea behind the definition \eqref{def:tau(u)} and
\thmref{thm:commutativity} is that the linear functional $\Phi(H;v)$ and
$(\hat\tau(u)\Phi)(H;v)$ are analogues of an $N$-point function of the
WZW model and the correlation funcion of the energy-momentum tensor,
respectively.

First we define the geometric data on which the WZW model lives. Let
$\X = S \times X_\tau$ be the trivial family of elliptic curves and
$\pi_{\X/S}:\X \to S$ be the projection. The divisor of $\X$
corresponding to the parameter $z_i$ is denoted by $P_i$ and their sum
by $D$:
\begin{equation}
    P_i := S \times z_i \bmod{\Integer + \Integer\tau} \subset \X, \qquad
    D := P_1 + \cdots + P_N
\label{def:divisor}
\end{equation}

The following general definitions of $N$-point functions is valid not
only for the WZW model but also for the free field theories which we use
later.

Assume that the following Lie algebraic data on $\X$ are given:
\begin{itemize}
\item $\fraka_\X$: a Lie algebra bundle with a fibre isomorphic to a Lie
      algebra $\fraka$.

\item $\langle \cdot \mid \cdot \rangle: \fraka_\X \times_{\O_S}
      \fraka_\X \to \Omega_{\X/S}^1$: an $\O_S$-bilinear
      $\Omega_{\X/S}^1$-valued pairing which is a 2-cocycle up to exact
      forms:
\begin{equation}
\begin{gathered}
    \fraka_\X \times \fraka_\X \owns (A, B) \mapsto
    \langle A \mid B \rangle \in \Omega_{\X/S}^1, \\
    \langle A \mid B \rangle + \langle B \mid A \rangle 
    \in d_{\X/S} \O_\X, \\
    \langle [A,B] \mid C \rangle + \langle [B,C] \mid A \rangle 
    + \langle [C,A] \mid B \rangle
    \in d_{\X/S} \O_\X,
\end{gathered}
\label{pairing:sheaf}
\end{equation}
     for any $A,B,C \in \fraka_\X$. For example, if there is a
     connection $\nabla$ on $\fraka_\X$ along the fibre compatible with
     the Lie algebra structure,
\begin{equation}
\begin{gathered}
    \nabla:\fraka_\X\to\fraka_\X\tensor_{\O_\X}\Omega_{\X/S}^1,\\
    \nabla[A,B] = [\nabla A, B] + [A, \nabla B] 
    \text{ for }A,B \in \fraka_\X,
\end{gathered}
\label{connection:sheaf}
\end{equation}
    and an invariant $\O_\X$-inner product $(\cdot \mid \cdot)$ of
    $\fraka_\X$ which satisfies
\begin{equation}
    d_{\X/S}(A|B) = (\nabla A \mid B) + (A\mid\nabla B) 
    \in \Omega_{\X/S}^1
    \quad
    \text{for $A,B\in\fraka_\X$},
\label{compati-conn-innerprod:sheaf}
\end{equation}
     then $\langle A \mid B \rangle := (\nabla A \mid B )$ has desired
     properties.
\end{itemize}

\begin{example}
\label{ex:wzw:bundle}
 In this section we use the following data. The Lie algebra bundle
 $\g_\X$ over $\X$ is defined as the quotient of $S \times \Complex
 \times \g$ by the $\Integer^2$-action,
\begin{equation}
    (m,n)\cdot (H;t;A)
    =
    (H;t+m\tau+n; e^{2\pi i m \ad H} A ),
\label{def:Z2-action:gX}
\end{equation}
 for $(m,n)\in \Integer^2$ and $(H;t;A) \in S \times \Complex \times
 \g$. The connection $\nabla$ of \eqref{connection:sheaf} is the trivial
 differentiation $d/dt \tensor dt$ along the elliptic curve and the
 invariant bilinear form $(\cdot\mid\cdot)$ is defined by
 \eqref{def:killing}.
\end{example}

For $i=1,\dots,N$, we put
\begin{equation}
\begin{split}
    \fraka_S^{P_i} 
    &:= (\pi_{\X/S})_\ast (\fraka_\X(\ast P_i))^{\wedge}_{P_i}
    \simeqq (\fraka\tensor \O_S)((x_i)),
\\
    \fraka_{S,+}^{P_i} 
    &:= (\pi_{\X/S})_\ast(\fraka_\X)^{\wedge}_{P_i}
    \simeqq (\fraka\tensor \O_S)[[x_i]],
\\
    \fraka_S^D &:=\bigoplus_{i=1}^N\fraka_S^{P_i}, \qquad
    \hat\fraka_S^D := \fraka_S^D \oplus \O_S \khat,
\end{split}
\label{def:loop-alg:sheaf}
\end{equation}
where $x_i$ is a formal local parameter and the central extension
$\hat\fraka_S^D$ is defined by the $\O_S$-valued 2-cocycle
\begin{equation}
    c(A,B) 
    := \sum_{i=1}^N \Res_{P_i} \langle A_i \mid B_i \rangle,
\label{def:cocycle:sheaf}
\end{equation}
where $A=(A_i)_{i=1}^N, B=(B_i)_{i=1}^N \in \fraka_S^D$

\begin{example}
\label{ex:wzw:affine}
 For \exref{ex:wzw:bundle}, the central extension $\hat\g_S^{P_i}$ is
 naturally isomorphic to the space of $\hat\g$-valued functions on $S$,
 $\hat\g_S = \hat\g \tensor_\Complex \O_S$, where $\hat\g$ is the affine
 Lie algebra associated to $\g$.
\end{example}

Let $\aDPr$ be the space of meromorphic sections of $\fraka_\X$ which
are globally defined along the fibre of $\pi_{\X/S}$ and holomorphic
except at $D$:
\begin{equation}
    \aDPr := (\pi_{\X/S})_\ast (\fraka_\X(\ast D)),
\label{def:aDPr}
\end{equation}
which is naturally regarded as a Lie subalgebra of $\hat\fraka_S^D$.

\begin{example}
\label{ex:wzw:global}
 For the Lie algebra bundle $\g_\X$ in \exref{ex:wzw:bundle}, we denote
 $\aDPr$ by $\g^D_\Prin$. The section $w_{\alpha(H)}(t-z_i)e_\alpha$
 (cf.\ \eqref{def:w,zeta}) belongs to $\g^D_\Prin$.
\end{example}

Let us take $\hat\fraka_S^{P_i}$-modules $\M_i$ ($i=1,\dots,N$) with the
same level $\khat = k$ and define $\M = \M_1 \tensor \cdots \tensor
\M_N$. 

\begin{defn}
\label{defn:N-point-func} 
(i)
 The sheaf of {\em conformal blocks} $\sCB(\fraka_\X,D,\M)$ is defined
 to be the space of $\O_S$-linear functionals on $\M$ which vanishes on
 $\aDPr \M$: $\Phi \in \HOM_{\O_S}(\M,\O_S)$ belongs to
 $\sCB(\fraka_\X,D,\M)$ if and only if it satisfies
\begin{equation}
    \Phi(A_\Prin v) = 0
    \quad\text{for all $A_\Prin \in \aDPr$ and $v \in \M$.}
\label{ward-cb:sheaf}
\end{equation}
 This equation \eqref{ward-cb:sheaf} is called the {\em Ward identity}.

(ii)
 There is a flat connection on $\sCB(\fraka_\X,D,\M)$. A flat section is
 called the {\em $N$-point function}. The set of $N$-point functions is
 denoted by $\sCB^\hor(\fraka_\X,D,\M)$.
 
\end{defn}

The flat connection on $\sCB(\g_\X,D,\M)$ is defined as follows. We use
the notation,
\begin{equation}
    \rho_i^\ast(h_r\{\theta_i(x_i)\})
    = \sum_{m\in\Integer} \theta_{i,m} \rho_i^\ast(h_r[m]),
\label{def:h(theta)}
\end{equation}
for an element
$
    \theta_i(x_i)
    = \sum_{m\in\Integer} \theta_{i,m} x_i^m
$,
of $\O_S((x_i))$. The connection on $\sCB(\g_\X,D,\M)$ in the direction
of $\xi_r$ is defined by
\begin{equation}
    \nabla^\ast_{\der/\der \xi_r}
    = \frac{\der}{\der \xi_r} - \rho^\ast(h_r\{Z(t)\})
    := \frac{\der}{\der \xi_r} 
     - \sum_{i=1}^N \rho^\ast_i(h_r\{Z(x_i+z_i)\}),
\label{flat-conn:wzw:cb:xi}
\end{equation}
where $Z(t)$ is a meromorphic function with poles in $D$ and has
quasi-periodicity 
\begin{equation}
    Z(t+m\tau + n) = Z(t) - 2 \pi i m.
\label{def:Z}
\end{equation}
(In \eqref{flat-conn:wzw:cb:xi} $Z(t)$ is expanded around $z_i$ in the
power series of $x_i = t-z_i$ and substituted into
\eqref{def:h(theta)}.) For example, we can take $Z(t) =
\zeta_{11}(t-z_1)$ (cf.\ \eqref{def:w,zeta}).

We denote $\sCB(\g_\X,D,\M)$ and $\sCB^\hor(\g_\X,D,\M)$ by $\sCB(D,\M)$
and $\sCB^\hor(D,\M)$ for short. We mean by the {\em WZW model} the
theory on $\sCB(D,\M)$ or $\sCB^\hor(D,\M)$. Hereafter we assume that
the $\hat\g$-module $M_i$ is generated over $\hat\g$ by a $\g$-submodule
$V_i$ on which the centre $\hat k\in\hat\g$ acts as multiplication by
$k$ and $A\tensor x^m$ ($A\in\g$, $m>0$) acts by $0$. Put $\M_i := M_i
\tensor \O_S$, $V = V_1 \tensor \cdots \tensor V_N$.

\begin{lem}
\label{lem:N-point-func|V} 
 An $N$-point function $\Phi(H;v)$ of the WZW model is determined by its
 values on $v\in V(0)$, where $V(0)$ is the $0$-weight space of $V$. In
 other words, the following restriction map is injective:
\begin{equation}
    \sCB^\hor(D,\M) \owns \Phi \mapsto 
    \Phi(H;v) \in V^\ast(0) \tensor \O_S.
\label{N-point-func|V}
\end{equation}

\end{lem}

This is a consequence of the Ward identity \eqref{ward-cb:sheaf} and the 
flatness condition. 

Let us return to the discussion for the general Lie algebra bundle
$\fraka_\X$. Let $Q$ be a point not contained in $D$ and $\Vac_{Q,k}$ be
the $\hat\fraka_S^Q$-module induced from the trivial
$\hat\fraka^Q_{S,+}$-module $\O_{S,k}$ of level $k$:
\begin{equation}
    \Vac_{Q,k} := \Ind_{\hat\fraka^Q_{S,+}}^{\hat\fraka^Q_S} \O_{S,k}
\label{def:vac}
\end{equation}
where $\O_{S,k} = \O_S |0\rangle$ as a linear space, $\fraka_{S,+}^Q
\O_{S,k} = 0$ and $\hat k$ acts as a multiplication by $k$. We call
$\Vac_{Q,k}$ the {\em vacuum module} of level $k$ at $Q$.

\begin{lem}
\label{lem:propagation-vac}
The canonical linear map $\M \owns v \mapsto v \tensor |0\rangle \in \M
\tensor \Vac_{Q,k}$ induces isomorphisms:
\begin{equation}
\begin{split}
    \sCB(\fraka_\X,D,\M) &\isoot 
    \sCB(\fraka_\X,D+Q, \M \tensor \Vac_{Q,k}),
\\
    \sCB^\hor(\fraka_\X,D,\M) &\isoot 
    \sCB^\hor(\fraka_\X,D+Q, \M \tensor \Vac_{Q,k}).
\end{split}
\label{propagation-vac}
\end{equation}
\end{lem}

This property is called {\em propagation of vacua} in \cite{t-u-y:89}.

The following proposition reveals the real nature of the operator
$\hat\tau(u)$ in \eqref{def:tau(u)}. The operator obtained directly from 
the WZW model differs from $\hat\tau(u)$ by conjugation. Let us define
an operator $\tilde\tau(u)$ by
\begin{equation}
\begin{split}
    \tilde\tau(u) &= \Pi(H,\tau) \hat\tau(u) \Pi(H,\tau)^{-1}
\\
    &=
    \frac{1}{2} \sum_{r=1}^l \nabla_{\xi_r,u}^2
    + \sum_{r=1}^l
    \frac{\partial}{\partial \xi_r}\log\Pi(H,\tau) \nabla_{\xi_r,u}
\\  &+
    \frac{1}{2} \sum_{i,j=1}^N \sum_{\alpha\in\Delta}
    w_{\alpha(H)}(z_i-u) w_{-\alpha(H)}(z_j-u) 
    \rho_j^\ast(e_{-\alpha}) \rho_i^\ast(e_{\alpha})
    + 2 \pi i \hvee \frac{\der}{\der \tau} \log \Pi(H,\tau),
\end{split}
\label{def:tilde-tau(u)}
\end{equation}
where $\Pi(H,\tau)$ is the normalised Weyl-Kac denominator,
\begin{equation}
    \Pi(H,\tau) =
    q^{\dim\g/24} (q;q)_\infty^l
    \prod_{\alpha\in\Delta_+} (e^{\pi i \alpha(H)}-e^{-i\pi\alpha(H)})
    \prod_{\alpha\in\Delta} (q e^{2\pi i\alpha(H)};q)_\infty.
\label{def:Pi}
\end{equation}
Here we use the usual notations, $q = e^{2\pi i \tau}$ and $(x;q)_\infty 
= \prod_{n=0}^\infty (1-xq^n)$. 

\begin{prop}
\label{prop:tau(u)=correl-func}
 According to \lemref{lem:propagation-vac}, there is a $(N+1)$-point
 function $\tilde\Phi$ corresponding to $\Phi \in \sCB^\hor(D,\M)$. We
 have 
\begin{equation}
    (\tilde\tau(u) \Phi) (H;v)
    =
    \tilde\Phi(H; v \tensor S[-2]|0\rangle),
\label{tau(u)=correl}
\end{equation}
 for $H \in S$, $v\in V(0)$, where $u$ is the coordinate of $Q$ on the
 complex plane and $S[-2]$ is defined as a coefficient of the {\em
 Sugawara tensor},
\begin{equation}
    S(u) := \frac{1}{2} \sum_{p=1}^{\dim\g} \NP J_p(u) J_p(u) \NP
    = \sum_{n\in\Integer} S[n] z^{-n-1}.
\label{def:sugawara}
\end{equation}
 Here $\{J_p\}_{p=1,\dots,\dim\g}$ is an orthonormal basis of $\g$ and
 the symbol $\NP\ \NP$ is the normal ordering operation.
\end{prop}

This means that $\tilde\tau(u)\Phi(v_1\tensor \cdots \tensor v_N)$ is
the correlation function $\langle S(u) v_1(z_1) \cdots v_N(z_n)\rangle$
of $S(u)$ in the context of the conformal field theory.

Hitherto the level is arbitrary. To prove \thmref{thm:commutativity}, we
need to fix the level to the critical value, $k = -\hvee$, where
$S[-2]|0\rangle$ is a singular vector of imaginary weight. Roughly
speaking, by virtue of this fact the correlation function of two Sugawara
tensors, $\langle S(u) S(u') v_1(z_1)\cdots v_N(z_n)\rangle$, is
irrelevant to the order of insertion of $S(u)$ and $S(u')$, from which
the commutativity \eqref{[tau,tau]=0} follows.

%%%%%%%%%%%%%%%%%%%%%%%%%%%%%%%%%%%%%%%%%%%%%%%%%%%%%%%%%%%%%%%%%%%%%%
\section{Wakimoto modules at the critical level}
\label{sec:wakimoto-critical}

The Bethe vector \eqref{BA-vector} is constructed by means of the
Wakimoto realisation of affine Lie algebras from the free field theory.
In this section we review basic facts about the Wakimoto
representations, following \cite{kur:91}. See also \cite{wak:86},
\cite{fei-fre:88}, \cite{fei-fre:90}, \cite{f-f-r:94}.

The {\em bosonic ghost fields},
\begin{equation}
    \beta_\alpha(z) = \sum_{m\in\Integer} z^{-m-1}\beta_\alpha[m], \qquad 
    \gamma^\alpha(z) = \sum_{m\in\Integer} z^{-m}\gamma^\alpha[m], \qquad 
    (\alpha \in \Delta_+)
\label{def:beta-gamma}
\end{equation}
satisfy the following operator product expansions:
\begin{equation}
    \beta_\alpha(z) \gamma^{\alpha'}(w) 
    \sim \frac{\delta_\alpha^{\alpha'}}{z-w}.
\label{ope:ghost}
\end{equation}
We denote the Heisenberg algebra generated by $\beta_\alpha[m]$ and
$\gamma^\alpha[m]$ ($\alpha\in\Delta_+$, $m\in\Integer$) by
$\widehat\Gh(\g)$.

The {\em ghost Fock space} $\ghFock$ is defined as a left
$\widehat\Gh(\g)$-module generated by the vacuum vector $\ghvac$,
satisfying
\begin{equation}
    \beta_\alpha[m]  \ghvac = 0, \qquad 
    \gamma^\alpha[n] \ghvac = 0
\label{def:ghost-vac}
\end{equation}
for any $\alpha\in\Delta_+$, $m\geqq 0$, $n>0$. The {\em normal ordered
product} $\np P \np$ of a monomial $P$ of $\beta_\alpha[m]$'s and
$\gamma^\alpha[m]$'s is defined by putting annihilation operators of
$\ghvac$ appearing in $P$ to the right side in the product.

The {\em free boson fields},
\begin{equation}
\begin{gathered}
    \phi_i(z) := \phi_i[0] \log z +
    \sum_{m\in\Integer\setminus\{0\}} \frac{z^{-m}}{-m} \phi_i[m], \qquad 
    \partial\phi_i(z) := \sum_{m\in\Integer} z^{-m-1} \phi_i[m], \\
    \phi(H;z) = \sum_{i=1}^l a_i \phi_i(z), \qquad \phi[H;m] =
    \sum_{i=1}^l a_i \phi_i[m],
\end{gathered}
\label{def:boson}
\end{equation}
for $H= \sum_{i=1}^l a_i H_i \in\h$ have trivial operator product
expansions: 
\begin{equation}
    \phi(H;z) \phi(H';w) \sim 0, \qquad
    \partial\phi(H;z) \partial\phi(H';w) \sim 0,
\label{ope:boson}
\end{equation}
for any $H, H' \in \h$. 

The commutative algebra generated by $\phi_i[m]$ ($i=1,\dots,l$,
$m\in\Integer$) is denoted by $\widehat\Bos(\g)$. 

For any one-form $\lambda(x)dx \in \h^\ast \tensor \Complex((x))
dx$, we define a one-dimensional representation
$\sigma_{\lambda(x)dx}$ of $\widehat\Bos(\g)$ by
\begin{equation}
\begin{gathered}
    \sigma_{\lambda(x)dx} = \Complex |\lambda(x)dx\rangle,
\\
    f(x)\, |\lambda(x)dx\rangle =
    \Res_{x=0} (\lambda(x), f(x)) \,dx,
\end{gathered}
\label{def:sigma-lambda}
\end{equation}
where $f(x) \in \h \tensor \Complex((x))$ is identified with an
element of $\widehat\Bos$ by the isomorphism defined by $H_i \tensor
x^m \mapsto H_i[m]$ and $(\cdot, \cdot)$ is the canonical pairing of
$\h^\ast$ and $\h$. In other words,
\begin{equation}
    \phi[H;m] \,|\lambda(x)dx\rangle 
    = \lambda^{(-m-1)}(H) \, |\lambda(x)dx\rangle,
\label{phi[H;m]|lambda>}
\end{equation}
where $\lambda(x)dx = \sum_{n\in\Integer} \lambda^{(n)} x^n \,dx$,
$\lambda^{(n)} \in \h^\ast$. Hereafter we assume that $\lambda^{(n)} =
0$ for $n \leqq -2$.

\begin{prop}{%
\cite{wak:86}, \cite{fei-fre:88}, \cite{fei-fre:90}, \cite{kur:91}.}
\label{prop:wakimoto}
 For each Chevalley generator $H_i$, $E_i$ or $F_i$ of $\g$, there
 exists a differential polynomial of the free fields,
\begin{equation*}
    X(z) =
    \sum_{m\in\Integer} X[m]z^{-m-1} := \np
     R(X;\gamma(z),\beta(z),\partial\phi(z)) \np,
\end{equation*}
 which gives the corresponding Kac-Moody current. Namely, a Lie algebra
 homomorphism $\omega$ from the affine Lie algebra
 $\hat\g=\g\tensor\Complex[t,t^{-1}]\oplus \Complex\hat k$ to a
 completion of $\widehat\Gh(\g) \tensor \widehat\Bos(\g)$ can be defined
 by
\begin{equation}
    \omega(X\tensor t^m) = X[m], \qquad
    \omega(\hat k) = - h^\vee,
\label{wakimoto-hom-km}
\end{equation}
 for all $X\in\g$, $m\in\Integer$, where $\hat k$ is the centre of
 $\hat\g$ and $h^\vee$ is the dual Coxeter number of $\g$. Moreover, the
 Sugawara tensor $S(z)$ defined by \eqref{def:sugawara}
 is expressed in terms of the free bosons as follows:
\begin{equation}
    S(z) :=
    \frac{1}{2}
    \sum_{r=1}^l \np \partial\phi(h_r;z) \partial\phi(h_r;z)\np 
    -
    \frac{1}{2} \partial^2\phi(2\rho;z),
\label{def:Sug(z)}
\end{equation}
\end{prop}

\begin{defn}
\label{def:wakimoto-rep:critical}
Denote $\ghFock \tensor \sigma_{\lambda(x)dx}$ by
 $\Waki_{\lambda(x)dx}$. We regard this module as a
 $\hat\g$-module through the above homomorphism $\omega: \hat\g \to
 \widehat\Gh(\g) \tensor \widehat\Bos(\g)$ and call it a {\em Wakimoto
 module} of level $-h^\vee$ (or of critical level) with highest weight
 $\lambda(x)dx$.
\end{defn}

The Wakimoto module contains a $\g$-submodule isomorphic to the dual
Verma module $M^\ast_\lambda$ of $\g$ with the highest weight
$\lambda=\lambda^{(-1)}$.  (See Proposition 4.4 of \cite{kur:91}.) We
denote it by $\Waki^{0}_{\lambda(x)dx}$.  It satisfies for any $m>0$ and
$X\in\g$,
\begin{equation}
    X[m] \Waki^{0}_{\lambda(x)dx} = 0.
\label{nilp=zero:wakimoto}
\end{equation}

%%%%%%%%%%%%%%%%%%%%%%%%%%%%%%%%%%%%%%%%%%%%%%%%%%%%%%%%%%%%%%%%%%%%%%
\section{Free field theories and Bethe vectors}
\label{sec:free-field,ba}

We can ``decompose'' the WZW model in \secref{sec:wzw-critical} into
free field theories, using the Wakimoto realisation,
\thmref{prop:wakimoto}. The Bethe vectors in \thmref{thm:BA} are nothing
but the $N$-point functions of the free field theories.

Let us define free field theories in the framework introduced in
\secref{sec:wzw-critical}. We defined the Lie algebra bundle $\g_\X$ as
a quotient of $S \times \Complex \times \g$ by the $\Integer^2$-action
\eqref{def:Z2-action:gX}. Note that this action preserves the triangular 
decomposition $\g = \n_+ \oplus \h \oplus \n_-$. We regard {\em vector}
bundles,
\begin{align}
    \bbeta_\X  &:= \Integer^2 \backslash S \times \Complex \times \n_+ 
\label{def:betaX}
\\
    \bgamma_\X &:= 
    (\Integer^2 \backslash S \times \Complex \times \n_-) 
    \tensor_{\O_\X}
    \Omega^1_{\X/S},
\label{def:gammaX}
\\
    \Bos_\X &:= \Integer^2 \backslash S \times \Complex \times \h,
\label{def:bosX}
\end{align}
as abelian Lie algebra bundles and put $\Gh_\X := \bbeta_\X \oplus
\bgamma_\X$. We apply \defref{defn:N-point-func} to $\fraka_\X = \Gh_\X$ 
and $\fraka_\X = \Bos_\X$. The pairing \eqref{pairing:sheaf} for
$\Bos_\X$ is trivial (i.e., $\langle \cdot \mid \cdot \rangle = 0$) and 
the pairing for $\Gh_\X$ is defined by
\begin{equation}
    \langle (A_1,B_1 dt) \mid (A_2,B_2 dt) \rangle^\gh
    :=
    ( A_1 \mid B_2 ) dt - ( B_1 \mid A_2 ) dt \in \Omega_\X^1,
\label{def:pairing:ghost}
\end{equation}
where $A_i \in \bbeta_\X$, $B_i dt \in \bgamma_\X$ for $i=1,2$ and
$(\cdot|\cdot)$ denotes the inner product defined by
\eqref{def:killing}. Here we identify $\bbeta_\X$ and $\bgamma_\X$
with a subbundle of $\g_\X$ and a subbundle of $\g_\X \tensor
\Omega^1_{\X/S}$ respectively.

The algebra defined by \eqref{def:loop-alg:sheaf} for $D=P$ (a point) in
this case is isomorphic to $\widehat\Gh_S(\g) =
\widehat\Gh(\g)\tensor_\Complex\O_S$ when $\fraka_\X = \Gh_\X$ and to
$\widehat\Bos_S(\g) = \widehat\Bos(\g)\tensor_\Complex\O_S$ when
$\fraka_\X = \Bos_\X$.

We denote the sheaf of conformal blocks $\sCB(\fraka_\X,D,\M)$ and the
space of $N$-point functions $\sCB^\hor(\fraka_\X,D,\M)$ defined in
\defref{defn:N-point-func} for $\fraka_\X = \Gh_\X, \Bos_\X$ by
\begin{align}
    &\sCB^{\gh}(D,\M) := \sCB(\Gh_\X,D,M),& 
    &\sCB^\bos (D,\M) := \sCB(\Bos_\X,D,M),
\label{def:gh/bos-cb}
\\
    &\sCB^{\gh,\hor}  (D,\M) := \sCB^\hor(\Gh_\X,D,M),& 
    &\sCB^{\bos,\hor} (D,\M) := \sCB^\hor(\Bos_\X,D,M),
\label{def:gh/bos-N-point}
\end{align}
and define
\begin{equation}
\begin{split}
    \sCB^\free(D,\M) &:= \sCB(\Gh_\X\oplus\Bos_\X,D,\M),
\\
    \sCB^{\free,\hor}(D,\M) &:= \sCB^\hor(\Gh_\X\oplus\Bos_\X,D,\M).
\end{split}
\label{def:free-cb/N-point}
\end{equation}

Assume that $\M_i$ is a $\widehat\Gh_S^{P_i}$-module and $\N_i$ is a
$\widehat\Bos_S^{P_i}$-module. Then $\M_i\tensor_{\O_S} \N_i$ is a
$(\widehat\Gh_S^{P_i}\oplus \widehat\Bos_S^{P_i})$-module and hence is a
$(\Gh_\X\oplus\Bos_\X)_S^{\wedge,P_i}$-module. The important fact is
that the $N$-point function of the free field theory naturally gives the
$N$-point function of the WZW model.
\begin{prop}
\label{prop:wzw-cc,cb<->free-cc,cb:sheaf}
 $\gDPr (\bigotimes(\M_i\tensor_{\O_S}\N_i)) \subset
 (\Gh^D_\Prin\oplus\Bos^D_\Prin)(\bigotimes(\M_i\tensor_{\O_S}\N_i))$. Hence
 the identity morphism from $\bigotimes (\M_i\tensor_{\O_S} \N_i)$ to
 itself induces an $\O_S$-linear map:
\begin{equation}
    \iota: \sCB^\free(D,\bigotimes (\M_i\tensor_{\O_S} \N_i)) \to
            \sCB(D,\bigotimes (\M_i\tensor_{\O_S} \N_i)).
\end{equation}
 Moreover, this induces a $\Complex$-linear map between spaces of
 $N$-point functions:
\begin{equation}
    \iota: \sCB^{\free,\hor}(D,\bigotimes (\M_i\tensor_{\O_S} \N_i)) \to
            \sCB^{\hor}(D,\bigotimes (\M_i\tensor_{\O_S} \N_i)),
\end{equation}
\end{prop}
If $\M_i = \ghFock \tensor \O_S$ and $\N_i = \sigma_{\mu_i(x)dx}
\tensor \O_S$ for $\mu_i(x) dx \in \h^\ast((x)) dx$, then an
$N$-point function of the free field gives an $N$-point function of the
WZW model with the Wakimoto modules. It is easy to see that
\begin{equation}
    \sCB^{\free,\hor}(D, \bigotimes (\M_i \tensor_{\O_S} \N_i))
    \simeqq
    \sCB^{\gh,\hor}(D, \M) \tensor_{\Complex} \sCB^{\bos,\hor}(D, \N),
\label{free=gh*bos}
\end{equation}
hence it is enough to find $N$-point functions in
$\sCB^{\gh,\hor}(D,(\ghFock)^{\tensor N})$ and
$\sCB^{\bos,\hor}(D,\bigotimes_{i=1}^N \sigma_{\mu_i(x)dx} \tensor
\O_S)$ to construct an $N$-point function of the WZW model. Let us
denote $\bsigma_{\vec \mu dx} := \bigotimes_{i=1}^N \sigma_{\mu_i(x)dx}
\tensor \O_S$ and $|\vec\mu dx\rangle := \bigotimes_{i=1}^N
|\mu_i(x)dx\rangle$ for simplicity.

\begin{lem}
\label{lem:cb:free-fields}
(i)
$
    \sCB^{\gh,\hor}(D,(\ghFock)^{\tensor N}) = \Complex \Phi^\gh(H;v),
$
where
\begin{equation*}
    \Phi^\gh(H;v) 
    = \Pi(H,\tau)^{-1} 
    (\text{\rm coefficient of }(\ghvac)^{\tensor N}\text{\rm\ in }v).
\end{equation*}

(ii)
 $\sCB^{\bos,\hor}(D,\bsigma_{\vec \mu dx})$ is one-dimensional if and
 only if there exists $\mu(t)dt \in
 (\pi_{\X/S})_\ast(\h^\ast\tensor\Omega^1_\X(\ast D))$ such that each
 $\mu_i(x_i) dx_i$ is a Laurent expansion of $\mu(t)dt$ at $t=z_i$ with
 respect to $x_i = t-z_i$. In this case,
 $\sCB^{\bos,\hor}(D,\bsigma_{\vec \mu dx}) = \Complex \Phi^\bos(v)$,
 where
\begin{equation*}
    \Phi^\bos(v) 
    = (\text{\rm coefficient of }|\vec\mu dx\rangle\text{\rm\ in }v).
\end{equation*}
 Otherwise, $\sCB^\bos(D,\bsigma_{\vec \mu dx}) = 0$.
\end{lem}

Assume $\lambda_i$ ($i=1,\dots,N$) satisfies the condition
\eqref{charge-preserving}, and put $\mu_i = \lambda_i$ ($i=1,\dots,N$),
$\mu_{N+j} = -\alpha_{i(j)}$
($j=1,\dots,M$). \lemref{lem:cb:free-fields} guarantees that
$\sCB^{\free,\hor}(D+D', \bsigma_{\vec\mu dx})$ is one-dimensional if we
define $D = P_1 + \cdots + P_N$ ($P_i$ has the coordinate $z_i$), $D' =
Q_1+ \cdots + Q_M$ ($Q_j$ has the coordinate $t_j$) and $\mu_i(x)dx$
($i=1,\dots,N+M$) as the Laurent expansion of $\mu(t)dt$ defined by
\begin{equation}
    \mu(t)dt = \sum_{i=1}^N \lambda_i \zeta_{11}(t-z_i)
             - \sum_{j=1}^M \alpha_{i(j)} \zeta_{11}(t-t_j).
\label{def:mu(t)dt}
\end{equation}
The basis of this one-dimensional space is 
$\Phi^\free(H;v) = \Phi^\gh(H;v) \Phi^\bos(v)$. 

We assume that the parameters $z_i$ ($i=1,\dots, N$) and $t_j$
($j=1,\dots,M$) satisy the Bethe Ansatz equations \eqref{BA-eq}. Then by
Lemma 2 of \cite{f-f-r:94} (or by Corollary 5.2 of \cite{kur:91}), the
{\em screening vector} $\scr_j$ (cf.\ (5.21) of \cite{kur:91}) in
$\Waki_{\mu_{N+j}(x)dx}$ is a singular vector of imaginary
weight. Thanks to this property and \eqref{nilp=zero:wakimoto}, the
linear functional defined by
\begin{equation}
    \Psi(H;v) 
    := 
    \Phi^\free (H; v \tensor \scr_1 \tensor \cdots \tensor \scr_M)
\label{def:Psi}
\end{equation}
for $v \in \bigotimes_{i=1}^N \Waki^0_{\mu_i(x)dx} \simeqq
\bigotimes_{i=1}^N M_{\lambda_i}^\ast$ has the same property as
\eqref{tau(u)=correl}:
\begin{equation}
    \tilde\Psi(H;v \tensor S[-2]|0\rangle) 
    = \tilde\tau(u) \Psi(H;v).
\label{tau(u)=correl(Psi)}
\end{equation}
On the other hand, we can compute the left hand side of
\eqref{tau(u)=correl(Psi)}, using the expression \eqref{def:Sug(z)} and
the Ward identity \eqref{ward-cb:sheaf} for the free boson. The result
is
\begin{equation}
    \tilde\Psi(H;v \tensor S[-2]|0\rangle) 
    = \tau_\Psi(u) \Psi(H;v),
\label{Psi=eigen}
\end{equation}
namely $\Psi(H;v)$ is the eigenvector of $\hat\tau(u)$ with the
eigenvalue $\tau_\Psi(u)$ defined by \eqref{BA-eigenval}. 

The explicit form of $\Psi(H;v)$, \eqref{BA-vector}, is derived by the
same argument as that of the rational case. See Lemma 3 of
\cite{f-f-r:94}.

%%%%%%%%%%%%%%%%%%%%%%%%%%%%%%%%%%%%%%%%%%%%%%%%%%%%%%%%%%%%%%%%%%%%%%
\appendix
\section{Elliptic functions}
\label{app:elliptic-functions}

We follow the notations of \cite{mum} for theta functions and denote the 
odd theta function by
\begin{equation}
    \theta_{11}(z,\tau) =
    \sum_{n\in\Integer} 
    \exp\left(
       \pi i \tau \left(n + \frac{1}{2}\right)^2
    + 2\pi i \left(z + \frac{1}{2}\right)\left(n + \frac{1}{2}\right)
    \right).
\label{def:theta11}
\end{equation}

We use a multiplicatively and additively quasi-periodic function,
\begin{equation}
    w_c(z) := \frac{\theta_{11}'(0) \theta_{11}(z-c)}
                   {\theta_{11}(z) \theta_{11}(-c)}, \qquad
    \zeta_{11}(z) := \frac{d}{dz} \log \theta_{11}(z).
\label{def:w,zeta}
\end{equation}

These functions are characterised by the properties
\begin{gather}
    w_c(z+1) = w_c(z), \qquad
    w_c(z+\tau) = e^{2\pi i c} w_c(z), \qquad
    w_c(z) \sim z^{-1} \text{ around }z=0.
\label{wc(z):property}
\\
    \zeta_{11}(z+1) = \zeta_{11}(z), \qquad
    \zeta_{11}(z+\tau) = \zeta_{11}(z) - 2\pi i, \qquad
    \zeta_{11}(z)\sim z^{-1} \text{ around }z=0.
\label{zeta(z):property}
\end{gather}

%%%%%%%%%%%%%%%%%%%%%%%%%%%%%%%%%%%%%%%%%%%%%%%%%%%%%%%%%%%%%%%%%%%%%%%%%%%

\end{document}